\documentclass{amsart}

\oddsidemargin .25in \theoremstyle{plain}

\newtheorem{thm}{Theorem}%[section]

\newtheorem{cor}[thm]{Corollary}

\newtheorem{rem}[thm]{Remark}

\newtheorem{quest}[thm]{Question}

\newtheoremstyle{definition}{7pt plus6.3pt minus6.3pt}{7pt plus3pt minus3pt}%
{\rm}{}{\bf}{}{0.75em}{\thmname{#1}\thmnumber{ #2}\thmnote{\sl\stdspace#3}}
\theoremstyle{definition}\newtheorem{example}[thm]{Example}
\newtheorem{exercise}[thm]{\small Exercise}

\newcommand{\bbr}{\begin{rem}\em} 
\newcommand{\eer}{\end{rem}}

\newcommand{\bex}{\begin{example}} 
\newcommand{\eex}{\end{example}}

\newcommand{\bhw}{\begin{exercise}\small} 
\newcommand{\ehw}{\end{exercise}}

\newcommand{\be}{\begin{enumerate}}
\newcommand{\ee}{\end{enumerate}}

\def\R{\hbox{$\mathbb R$} }

\def\co{\colon\thinspace}

\def\dfn#1{{\em #1}}

\begin{document}

%%%%%%%%%%%%%%%%%%%%%%%%%%%%%%%%%%%%%%%%%%%%%%%%%%%%%%%%%%
\title{Contact structures on 3--manifolds are deformations of foliations}

\author{John B. Etnyre}
%\address{Stanford University, Stanford, CA 94305}
%\address{University of Pennsylvania, Philadelphia, PA 19104}
%\email{etnyre@math.upenn.edu}
%\urladdr{http://www.math.upenn.edu/\char126 etnyre}
\address{Georgia Institute of Technology, Atlanta, GA 30332-0160}
\email{etnyre@math.gatech.edu}
\urladdr{http://www.math.gatech.edu/\char126 etnyre}
\thanks{The author thanks Baptiste Chantraine  and the referee for helpful comments on the
first version of the paper. Supported in part by NSF CAREER Grant (DMS--0239600) and FRG-0244663.}

\begin{abstract}
In this note we observe, answering a question of Eliashberg and Thurston in \cite{EliashbergThurston98}, that all contact structures on a 
closed oriented 3--manifold are $C^\infty$-deformations of foliations.
\end{abstract}
%\keywords{tight, contact structure, Legendrian}
%\subjclass{Primary 53C15; Secondary 57M50}

\maketitle

%\tableofcontents

%%%%%%%%%%%%%%%%%%%%%%%%%%%%%%%%%%%%%%%%%%%%%%%%%%%%%%%%%%
\section{Introduction}
%%%%%%%%%%%%%%%%%%%%%%%%%%%%%%%%%%%%%%%%%%%%%%%%%%%%%%%%%% 
In \cite{EliashbergThurston98}, Eliashberg and Thurston proved that foliations could be approximated by contact structures. More precisely they 
established that % \begin{thm}[Eliashberg and Thurston, 1998 \cite{EliashbergThurston98}]
any co-dimension two $C^2$-foliation on a closed oriented 3--manifold $M,$ other than the product foliation of $S^2\times S^1,$ 
can be $C^0$-approximated by a positive and negative contact structure. 
%\end{thm}
This result brings up the natural question: Is every contact structure on a closed 3--manifold ``close'' to a foliation? To make sense of ``close'' we could
instead ask, as Eliashberg and Thurston did in \cite{EliashbergThurston98}: Is every contact structure on a closed 3--manifold a deformation of a foliation?
We say a contact structure $\zeta$ is the deformation of a foliation $\xi$ if there is a 1-parameter family of plane fields $\xi_t$ such that $\xi_0=\xi, \xi_1=\zeta$
and $\xi_t$ is a contact structure for $t>0.$ In this note we show that answer to this last question is indeed yes.
\begin{thm}\label{main}
Every positive and negative contact structure on a closed oriented 3--manifold 
is a $C^\infty$-deformation of a $C^\infty$-foliation. 
\end{thm}
The proof of this relies on the connection between open books and contact structures discovered by Giroux \cite{Giroux02}. Starting form open books
it is not hard to prove this theorem, but it does illuminate the nature of confoliations and the ``boundary'' of the space of contact structures in the space
of plane fields. In the last section we make a few comments concerning the theorem and its proof. 

%%%%%%%%%%%%%%%%%%%%%%%%%%%%%%%%%%%%%%%%%%%%%%%%%%%%%%%%%%
\section{Open books and contact structures}
%%%%%%%%%%%%%%%%%%%%%%%%%%%%%%%%%%%%%%%%%%%%%%%%%%%%%%%%%% 

Recall an open book decomposition of a closed oriented 3--manifold $M$ is a pair $(L,\pi)$ where $L$ is an oriented  
link in $M$ and $\pi\co (M\setminus L)\to S^1$
is a fibration with $\partial (\overline{\pi^{-1}(\theta)})=L$ for all $\theta\in S^1.$ We call $L$ the binding of the open book and the fibers of $\pi$ the
pages of the open book. Given an open book we can describe $M\setminus L$ as the mapping cylinder of a diffeomorphism $\psi\co \Sigma\to \Sigma$
of a surface $\Sigma.$ Indeed, we can recover $M$ and the open book $(L,\pi),$ up do diffeomorphism, from this data $(\Sigma,\psi).$ We call $\psi$ the
monodromy of the open book. For more on open books see \cite{EtnyreOBN}.

A contact structure $\xi$ is said to be supported by an open book $(L,\pi)$ if there is a contact 1-form $\alpha$ for $\xi$ such that
$\alpha(L)>0$ and $d\alpha$ is a volume form when restricted to each page of the open book. Thurston and Winkelnkemper, in \cite{ThurstonWinkelnkemper75},
showed that any open book supports some contact structure and it was observed by Giroux that this contact structure is unique. Moreover, Giroux
proved the following correspondence.
\begin{thm}[Giroux, 2002 \cite{Giroux02}]
If $M$ is a closed oriented 3--manifold then there is a one-to-one correspondence between oriented contact structures on $M$ up to isotopy and
open book decompositions of $M$ up to positive stabilization.
\end{thm}
It is not important for our arguments what the definition of positive stabilization is, so we do not define it here. Our main use of this theorem will be
the facts that (1) all contact structures are supported by open books and (2) the supported contact structure is unique up to isotopy. 

%%%%%%%%%%%%%%%%%%%%%%%%%%%%%%%%%%%%%%%%%%%%%%%%%%%%%%%%%%
\section{Proof of Theorem~\ref{main}}\label{construction}
%%%%%%%%%%%%%%%%%%%%%%%%%%%%%%%%%%%%%%%%%%%%%%%%%%%%%%%%%% 

To prove Theorem~\ref{main} we start with a positive contact structure $\zeta$ on a closed oriented 3--manifold $M$ and then choose some open book
$(L,\pi)$ for $M$ that supports $\zeta.$ (The proof for negative contact structures is similar. The details are left as an easy exercise for the reader.) 
We will then construct a foliation on $M$ associated to this open book and show that it can be perturbed into
a contact structure that is also supported by $(L,\pi). $ Thus the perturbed contact structures will have to be isotopic to $\zeta,$ confirming the 
theorem. We will denote the page and monodromy of the open book by $(\Sigma, \pi).$

One can construct an obvious foliation with Reeb components by replacing neighborhoods of the binding with Reeb components and
then ``spinning'' the pages of the open book so they limit to the Reeb components. We describe the procedure more carefully. 
Let $N$ be a neighborhood of one component of the binding. Choose coordinates $(r,\theta, \phi),$ so that the pages of the
open book intersect $N$ correspond to constant $\theta$ annuli and the binding corresponds to $r=0.$ Moreover assume
$N=\{(r,\theta, \phi) | r\leq 1+ 2\epsilon\}$ for some small fixed $\epsilon.$ We now choose two functions $\lambda(r)$ and $\delta(r)$ so
that $\lambda$ is zero on $[0,\frac 13],$ one for $r\geq 1,$ and strictly increasing on $[\frac 13, 1],$ and $\delta$ is zero on $[0,1],$ 
one for $r\geq 1+\epsilon$ and strictly increasing on $[1, 1+\epsilon].$ Now set 
\[
\alpha= \begin{cases}
\lambda(r)\, dr + (1-\lambda(r))\, d\phi  \quad &\text{for } r\leq 1,\\
\delta(r)\, d\theta + (1-\delta(r))\, dr &\text{for } r>1.
\end{cases}
\]
One may easily check that $\alpha\wedge d\alpha=0$ so $\xi=\ker\alpha$ gives a foliation of $N.$ The subset $N_1=\{(r,\theta,\phi)|r\leq 1\}$ of $N$ is
a Reeb component. (We will denote by $N_a$ the set $\{(r,\theta,\phi)| r\leq a\}.$) %See Figure~\ref{}. 
Note we can choose $\lambda(r)$ and $\delta(r)$ so that $\alpha$ defines a $C^\infty$-foliation on $N$.
It is clear that in $[1+\epsilon, 1+2\epsilon]\times T^2$ the foliation is by constant $\theta$ annuli. Thus
this foliation can be extended by the pages of the open book on $M\setminus N$ to a foliation of all of $M.$ (Of course, we insert this
standard model with a
Reeb component for each component of the binding $L$.) If we set $dz$ to the pull back of the coordinate $\vartheta$ on $S^1$ by the fibration
$\pi\co M\setminus L \to S^1$ then we can extend $\alpha$ by $dz$ to get a 1-form defining our foliation on all of $M.$  

We now show how to perturb this foliation in to a contact structure supported by the open book $(L,\pi).$ Again, we start by concentrating
on the neighborhood $N$ of a component of the binding $L.$ Here we set 
\[
\alpha_t= \alpha + t (r^2\, d\theta + (1+f(r))\, d\phi),
\]
where $f\co N\to \R$ is a function that is strictly decreasing, $f(0)=0, f(r)>-1$ for all $r$ and $f(r)<-1+\iota$ for all $r>1$ and $\iota$ some small number. 
Now we compute
\[
d\alpha_t=\begin{cases}
(tf'(r)-\lambda'(r))\, dr\wedge d\phi + t2r\, dr\wedge d\theta \quad&\text{for } r\leq 1,\\
(t2r+\delta'(r))\, dr\wedge d\theta + tf'(r)\, dr\wedge d\phi &\text{for } r>1.
\end{cases}
\]
Thus we have
\begin{align*}
&\alpha\wedge d\alpha_t =\\
& \begin{cases}
tr\bigl( 2[(1-\lambda(r))+t(1+f(r))]   -r(tf'(r) - \lambda'(r))\bigr)\, dr\wedge d\theta\wedge d\phi \quad &\text{for }  r\leq 1,\\
t\bigl( -f'(r)(\delta(r) + tr^2) + (1+f(r))(t2r +\delta'(r))  \bigr)\, dr\wedge d\theta\wedge d\phi & \text{for }  r> 1.
\end{cases}
\end{align*}
From this formula, it may easily be checked that $\alpha_t$ is a contact from on $N$ for $t>0.$ 

We would like to 
patch this into a family of 1-forms on $M\setminus N_{1+\epsilon}.$ To this end we briefly recall the Thurston-Winkelnkemper 
construction of contact forms on open books, \cite{ThurstonWinkelnkemper75}. We begin by thinking of $M\setminus N_{1+\epsilon}$ as the mapping cylinder  
of $\psi\co \Sigma \to \Sigma,$ the monodromy of the open book. That is 
\[
M\setminus N' = \Sigma\times[0,1]/\sim,
\]
where $(\psi(x), 0)\sim(x,1),$ and the coordinate on the $[0,1]$ factor is $z.$ 
We then find a 1-parameter family of 1-forms $\lambda_z$ on $\Sigma$ so that $d\lambda_z$ is a volume form on $\Sigma$ for all $z$ and
each $\lambda_z=(1+\epsilon+s)\, d\theta$ near each boundary component of $\Sigma,$ where we use ``polar'' coordinates $(s,\theta)$ near the boundary
component and the boundary corresponds to $s=0$ and $s$ is increasing into $\Sigma.$  Moreover, the $\lambda_z$ are chosen so
that they descend to give a form on $M\setminus N_{1+\epsilon}.$ Note the 1-form $dz,$ from the second paragraph of this section, corresponds
to $dz$ in these coordinates. The 1-form $\beta_t=dz + t\lambda_z$ will be a contact 1-form on $M\setminus N_{1+\epsilon}$ for small $t>0.$

To patch the two 1-forms together we consider the region $A=\overline{N\setminus N_{1+\epsilon}}.$ 
We use the above coordinates on $N$ as coordinates on $A.$ 
Near the boundary of $M\setminus N$  in $A$ the contact 1-form is 
$\beta_t=dz + t(1+\epsilon+s)\, d\theta.$ We use the map $\Psi(r, \theta, \phi)=(r-1-\epsilon, -\phi, \theta)$ to map $A\subset N$
to a neighborhood of the boundary of $M\setminus N_{1+\epsilon}.$ It is easy to check that this map is orientation preserving and when 
$N$ is glued to $M\setminus N_{1+\epsilon}$ 
using this map we recover $M.$ Pulling $\beta_t$ back to $A$ using this map we get $\Psi^*\beta_t= -tr\, d\phi + d\theta.$ We think of this form as defined 
only near 
$T_{1+2\epsilon}=\partial N_{1+2\epsilon}$ in $A.$ Similarly $\alpha_t=(1+tr^2)\, d\theta + t(1+f(r))\, d\phi$ is a form defined near $T_{1+\epsilon}=
\partial N_{1+\epsilon}$ in $A.$ In order to interpolate between these two forms we consider forms on $A$ of the type
\[
\gamma=g(r)\, d\phi + h(r)\, d\theta.
\]
This will be a contact form if and only if $g(r)h'(r)-h(r)g'(r)\not =0,$ that is if we think of $(g(r),f(r))$ as parameterizing a curve in $\R^2$ then the position vector
and velocity vector can never be co-linear. If we take $g(r)$ and $h(r)$ to be defined by $\Psi^*\beta_t$ and $\alpha_t$ near the boundary of $A$ then
we can clearly extend $g(r)$ and $h(r)$ to all of $A$ so that we have a contact form on $A.$ Moreover, it is easy to check that we can choose $g(r)$ so that
$g'(r)<0$ in $A.$

Let $\alpha_t$ be the 1-from on $M$ that equals $\alpha_t$ on $N_{1+\epsilon},$ $\beta_t$ on $M\setminus N$ and the 
from $g(r)\, dz + h(r)\, d\theta$ just constructed on
$A.$ This gives a well defined form for all $t\geq 0.$ Moreover, $\alpha_0$ is the form $\alpha$ above that defines the foliation $\xi$ and for small $t>0,$
$\alpha_t$ defines a contact structure $\xi_t=\ker \alpha_t.$ Thus the contact structure $\xi_t$ is clearly a deformation of the foliation $\xi.$ 

We are left to show that $\xi_t$ is supported by the open book $(L,\pi).$ For this we must see that $\alpha_t(L)>0$ and $d\alpha_t|_\text{page}$ is a 
volume form on $\Sigma.$ A component of $L$ corresponds to $r=0$ in $N$ and its positively oriented tangent vector is given by $\frac{\partial}{\partial \phi}.$
So $\alpha_t(\frac{\partial}{\partial \phi})=d\phi (\frac{\partial}{\partial \phi})=1>0.$ To check the second condition we consider the four regions 
$N_1, N_{1+\epsilon}\setminus N_1, A$ and
$M\setminus N.$ On $N_1$ the pages of the open book correspond to constant $\theta$ annuli. The form $d\alpha_t$ restricted to this annulus is
$(tf'(r) - \lambda'(r))\, dr\wedge d\phi.$ This form is never zero and the coefficient is always negative in $N_1,$ but the orientation on the annulus that
allows for $L$ to be properly oriented corresponds to the form $d\phi \wedge dr.$ So $d\alpha_t$ is a properly oriented non-zero 2-form on the pages in 
$N_1.$ Now on $N_{1+\epsilon}\setminus N_1$ the pages are still constant $\theta$ annuli and the 1-form restricts to $tf'(r)\, dr\wedge d\phi$ on these. Thus 
$\alpha_t$ is compatible withe the pages in this region. 
On $A$ the pages are again constant $\theta$ annuli, so the form restricted to this is $g'(r)\, dr\wedge d\phi.$ By the choice of $g$ above this is a properly 
oriented non-zero 2-form on the pages in $A.$ Finally in $M\setminus N$ it is clear $d\beta_t$ is a properly oriented non-zero 2-form on the pages by
construction (see \cite{ThurstonWinkelnkemper75}).

%%%%%%%%%%%%%%%%%%%%%%%%%%%%%%%%%%%%%%%%%%%%%%%%%%%%%%%%%%
\section{Comments}
%%%%%%%%%%%%%%%%%%%%%%%%%%%%%%%%%%%%%%%%%%%%%%%%%%%%%%%%%% 
In \cite{EliashbergThurston98} it was shown that any confoliation is isotopic through confoliations to an overtwisted contact structure. Using Theorem~\ref{main}
we can see other interesting facts concerning isotopies of confoliations. 
%Say a confoliation $\xi$ is of class $\mathcal{T},$ or is in $\mathcal{T},$ 
%if its Euler class $e(\xi)$ satisfies
%\[|\langle e(\xi), [\Sigma]\rangle|\leq -\chi(\Sigma),\]
%for a closed embedded surface $\Sigma\not=S^2,$
%\[|\langle e(\xi), [\Sigma]\rangle|=0,\]
%if $\Sigma=S^2$ and 
%\[\langle e(\xi), [\Sigma]\rangle\leq -\chi(\Sigma),\]
%if $\Sigma$ is a surface with boundary positively transverse to $\xi,$ and here the left hand side of the inequality refers to the Euler class of $\xi$ restricted to
%$\Sigma$ relative to any non-zero section of $T\Sigma\cap \xi$ along $\partial \Sigma.$
Eliashberg and Thurston generalized the notion of tightness to confoliations. They say a confoliation $\xi$ is \dfn{tight} 
if for every disk $D$ with $\partial D$
tangent to $\xi,$ but $D$ itself transverse to $\xi$ near $\partial D,$ there is a disk $D'$ with the same boundary as $D$ that 
is tangent to $\xi$ everywhere and satisfies $\langle e(\xi), [D\cup D']\rangle=0,$ where $e(\xi)$ is the Euler class of the plane field $\xi.$ 
If $\xi$ is a contact structure this agrees with the usual definition of
tight. If $\xi$ is a foliation the is equivalent to the foliation being Reebless. 
%If $\xi$ is a tight foliation of contact structure then it is in $\mathcal{T},$
%but it is not known if any tight confoliation is in $\mathcal{T}.$
Eliashberg and Thurston claimed that the perturbation of a tight foliation
is a tight contact structure (and in fact tight when pulled back to the universal cover). While there seems to be a gap in their argument, see \cite{Colin02}, the result is true, see
\cite{EtnyreLF2C}. 
 We have the immediate corollary of Theorem~\ref{main}.
\begin{cor}
If $\xi$ is a tight foliation  then we can isotop $\xi$ through tight confoliations to a foliation with Reeb components. 
In particular, a taut foliation is isotopic to a foliation with Reeb components via an isotopy through tight confoliations. 
\end{cor}

In the definition of the foliation in Section~\ref{construction} the reader may verify that the leaves of the foliations coming from the pages of the open book spiraled 
towards the Reeb component in a clockwise manner (as seen in a constant $\phi$ slice). We could have considered a foliation that spiraled in the opposite
direction or changed the direction of the spiraling in the Reeb component itself. Making either of these changes (but not both) would have resulted in a foliation
that perturbs to an overtwisted contact structure. (To see this take a meridional disk 
in $N$ of radius slightly larger than one. This will give a transverse unknot in the contact structure that violates the Bennequin inequality.) 
So if the open book supported a tight contact structure then this perturbed contact structure could not be supported by the open book. More generally, it is not too hard
to see that the perturbed contact structure will never be supported by the given open book.

Finally we observe that if a contact structure is virtually overtwisted ({\em i.e.} becomes overtwisted when pulled back to a finite cover) then any foliation of which it
is a perturbation must contain Reeb components, since Reeblees foliations are tight and perturb to tight contact structures.
Thus the foliations constructed in the proof Theorem~\ref{main} are
as good as can be expected for most contact structures. However, if a contact structure is universally tight ({\em i.e.} its pull back to
the universal cover is tight) then one may hope for better. In particular, we end with the following questions:
\begin{quest}
If $M$ has infinite fundamental group then 
is every universally tight contact structure the deformation of a Reebless foliation?
\end{quest}
\begin{quest}
If $M$ has infinite fundamental group then 
is every universally tight contact structure the deformation of a taut foliation?
\end{quest}
The condition on the fundamental group is necessary since a theorem of Novikov says that if the fundamental group is finite then there must be a Reeb component 
in the foliation. 
%\bibliography{references}
%%\bibliographystyle{gtart}
%\bibliographystyle{plain}
\def\cprime{$'$} \def\cprime{$'$}

\end{document}